\documentclass[a4paper,pdftex,Skim]{amsart}

\usepackage{geometry}
\usepackage{afterpage}
\usepackage{eurosym}

\usepackage{graphicx}
\graphicspath{{Figures/}}

\usepackage{courier}

\usepackage{hyperref}

\theoremstyle{plain}
\newtheorem{Theorem}{Theorem}[section]

\theoremstyle{definition}
\newtheorem{Definition}[Theorem]{Definition}

\theoremstyle{remark}

\newcommand{\TDI}[2]{\ensuremath{\operatorname{TDI}_{#1}(#2)}}
\newcommand{\TDC}[2]{\ensuremath{\operatorname{W}_{#1}(#2)}}
\newcommand{\FDC}[2]{\ensuremath{\operatorname{L}_{#1}(#2)}}
\begin{document}

\title[The Top-Dog index]{The Top-Dog Index: A New Measurement for the
  Demand Consistency of the Size Distribution in Pre-Pack Orders for a
  Fashion Discounter with Many Small Branches}

\begin{abstract}
  We propose the new Top-Dog-Index, a measure for the branch-dependent
  historic deviation of the supply data of apparel sizes from the
  sales data of a fashion discounter. A common approach is to
  estimate demand for sizes directly from the sales data. This approach may
  yield information for the demand for sizes if
  aggregated over all branches and products. However, as we will show in a
  real-world business case, this direct approach is in general
  not capable to provide information about each branch's individual
  demand for sizes: the supply per branch is so small that either the
  number of sales is statistically too small for a good estimate
  (early measurement) or there will be too much unsatisfied demand
  neglected in the sales data (late measurement). Moreover, in
  our real-world data we could not verify any of the demand
  distribution assumptions suggested in the literature.  Our approach
  cannot estimate the demand for sizes directly.  It can, however,
  individually measure for each branch the scarcest and the amplest sizes,
  aggregated over all products. This measurement can
  iteratively be used to adapt the size distributions in the pre-pack
  orders for the future.  A real-world blind study shows the potential
  of this distribution free heuristic optimization approach: The gross
  yield measured in percent of gross value was almost one percentage
  point higher in the test-group branches than in the control-group
  branches.
\end{abstract}

\author{Sascha Kurz}
\author{J\"org Rambau}
\author{J\"org Schl\"uchtermann}
\author{Rainer Wolf}

\address{Sascha Kurz\\Fakult\"at f\"ur Mathematik, Physik und Informatik\\Universit\"at Bayreuth\\Germany}
\email{sascha.kurz@uni-bayreuth.de}
\address{J\"org Rambau\\Fakult\"at f\"ur Mathematik, Physik und Informatik\\Universit\"at Bayreuth\\Germany}
\email{joerg.rambau@uni-bayreuth.de}
\address{J\"org Schl\"uchtermann\\Fakult\"at f\"ur Rechts- und Wirtschaftswissenschaften\\Universit\"at Bayreuth\\Germany}
\email{j.schluechtermann@uni-bayreuth.de}
\address{Rainer Wolf\\Betriebswirtschaftliches Forschungszentrum f\"ur Fragen der
  mittelst\"andischen Wirtschaft e.\,V.{} an der Universit\"at Bayreuth, Germany}
\email{rainer.wolf@uni-bayreuth.de}



\thanks{This research was
    supported by the \emph{Bayerische Forschungsstiftung, Project
      \emph{DISPO---A Decision Support for the Integrated Size and Price
        Optimization}}}

\maketitle

\section{Introduction}
\label{sec:introduction}

%% The economic success of a fashion discounter lies mostly in its ability to
%% predict the customers' demand for individual products.  
The financial performance of a fashion discounter depends very much on
its ability to predict the customers' demand for individual products.
More specifically: trade exactly what you can sell to your customers.
This task has two aspects: offer what your customers would like to
wear because the product as such is attractive to them and offer what
your customers can wear because it has the right size.

In this paper, we deal with the second aspect only: meet the demand for sizes
as accurately as possible.  The first aspect, demand for products, is a very
delicate issue: Products in a fashion discounter are never replenished because
of lead times of around three months.  Therefore, there will never be historic
sales data of an item at the time when the order has to be submitted (except
for the very few ``never-out-of-stock''-items, NOS items, for short).

When one considers the knowledge and experience of the professional buyers
employed at a fashion discounter---acquired by visiting expositions, reading
trade journals, and the like---it seems hard to imagine that a forecast for the
demand for a product could be implemented in an automated decision support
system at all. We seriously doubt that the success of fashion product can be
assessed by looking at historic sales data only. In contrast to this, the demand
for sizes may stay reasonably stable over time to extract useful information from
historic sales data.

In the historic sales data the influences of demand for products and demand
for sizes obviously interfere.  Moreover, it was observed at our partners'
branches that the demand for sizes seems not to be constant over all around
1\,200~branches.

The main question of this work is: how can we forecast the demand for sizes
individually for each branch or for a class of branches?

\subsection{Related Work}
\label{sec:introduction:related-work}

Interestingly enough, we have not found much work that exactly deals with our
task. It seems that, at first glance, the problem of determining the size
distribution in delivery pre-packs can be considered as simple regression once
you have historic sales data: Just estimate the historic size profile and fit
your delivery to that. At least two trivial US-patents
\cite{Rose+Leven+Woo:DemandForecasting:2006,Rose+Leven+Woo:AssortmentDecisions:2006}
have recently been granted and published along these lines (which witnesses that
the US patent system may not have employed the necessary expertise in their
patent evaluations \ldots).

In our problem, however, the historic sales data is \emph{not} necessarily
equal to the historic demand data, and it is interesting how to find the
demand data in the sales data in the presence of unsatisfied demand and very
small delivery volumes per branch and per product.

The type of research closest to ours seems to be classified as assortment
optimization.  In a sense, we want to decide on the start inventory level of
sizes in a pre-pack for an individual product in an individual branch.  (Let
us ignore for a moment that these unaggregated inventory levels are very small
compared to other inventory levels, e.g., for grocery items.)

Mostly, the successful approaches deal rather with NOS items than with
perishable and not replenishable fashion goods.  For example,
assortment optimization in the grocery sector
\cite[Section~4]{Koek+Fisher:AssortmentOptimization:2006}\,---\,one of
the very few papers documenting a field study\,---\,can usually
neglect the effects of stockout substitution in sales data, which make
demand estimation from sales data much more reliable.  There is work
on the specific influence of substitution on the \emph{optimization of
  expected profit} (see, e.g.,
\cite{Mahajan+vanRyzin:AssortmentsSubstitution:2001}), but the problem
of how to \emph{estimate demand parameters} from low-volume sales data
in the presence of stock-out substitution remains.

Much more work has been published in the field of dynamic pricing, where in
one line of research pricing and inventory decisions are linked.  See
\cite{Elmaghraby+Keskinocak:DynamicPricing-Inventory:2003,Chan+Shen+Simchi-Levi+Swann:PricingInventorySurvey:2006}
for surveys.

A common aspect of all cited papers (and papers cited there) that separates
their research from ours is the following: those papers, in some sense,
postulate the possibility to estimate a product's demand in an individual branch
directly from sales data, in particular from sales rates. In our real-world
application we have no replenishment, small delivery volumes per branch, lost
sales with unknown or even no substitution, and sales rates depending much more on
the success of the individual product at the time it was offered than on the
size. Therefore, estimating future absolute demand data from historic sales data
directly seems to need extra ideas, except maybe for the data aggregated over
many branches.

Size optimization can be found in a few offers of retail optimization systems,
like 7thOnline (\url{http://www.7thonline.com}).  It is not clear, on what kind
of research these products are based and under which assumptions they work
well.  Our partner firm has checked several offers in the past and did not
find any optimization tools that met their needs.

\subsection{Our contribution}
\label{sec:introduction:contribution}

The main result of this work is: a useful forecast for the demand for
sizes on the level of individual branches is too much to be asked for,
but historic information about which sizes have been the scarcest and
which sizes have been the amplest ones can be obtained by measuring
the new \emph{Top-Dog-Index (TDI)}.  The TDI can be utilized in a
dynamic heuristic optimization procedure, that adjusts the size
distributions in the branches' corresponding pre-packs accordingly
until the difference between the scarcest and the amplest sizes can
not be improved anymore.  The main benefit of the TDI: it measures the
consistency of historic supply with sizes with the historic demand for
sizes in a way that is not influenced by the attractivity of the
product itself.  This way, we can aggregate data over all products of
a product group, thereby curing the problem with small delivery
volumes per branch and product and size.

The potential of our TDI-approach is shown in a blind-study with
20~branches and one product group (womens' outer garments). Ten
branches randomly chosen from the 20~branches (test-branches) received
size pre-packs according to our heuristics' recommendations, ten
received unchanged supply (control-branches).  The result: One
percentage point increase of gross yield per merchandise value for the
test-branches against the control-branches.  A conservative
extrapolation of this result for our partner would already mean a
significant increase of gross yield.

We have not seen any field study of this type documented in the literature so
far. The only documented yield management studies in the apparel retail business
(e.g., for dynamic pricing policies) try to prove the success of their methods by
showing that they would have gained something on the set of data that was used to
estimate the parameters of the model
\cite{Bitran+Caldentey+Mondschein:ClearanceMarkdown:1998}. Such tests are very
far from reality: a by-product of our study is that in our case it makes
absolutely no sense to compare gross yield data of a fashion discounter across
seasons, because the differences between the yields in different seasons are
much larger than the differences caused by anything we are interested in. This
was the reason for us to use the parallel blind testing instead.

Further tests are planned by our partner on a larger scale in the near future.

\subsection{Outline of the paper}
\label{sec:introduction:outline}

In Section~\ref{sec:real-world-problem}, we briefly restate the
real-world problem we are concerned with.  In
Section~\ref{sec:data-analysis}, we show our experience with
straight-forward estimators for the demand distribution on sizes.  In
Section~\ref{sec:TDI}, we introduce the new Top Dog Index, which is
utilized in Section~\ref{sec:heuristic-optimization} in a heuristic
optimization procedure.  Section~\ref{sec:blind-study} is the
documentation of a field study containing a blind testing procedure
among two groups of branches: one supplied with and one supplied
without the suggestions from the first step of the optimization
heuristics.  We summarize the findings in
Section~\ref{sec:conclusion-and-outlook}, including some ideas for
further research.

\section{The real-world problem}
\label{sec:real-world-problem}

In this section, we state the problem we are concerned with.  Before
that we briefly provide the context in which our problem is embedded.

\subsection{The supply chain of a fashion discounter}
\label{sec:real-world-problem:supply-chain}

As in most other industries the overall philosophy of supply chain
management in fashion retailing is to coordinate the material flow
according to the market demand. The customer has to become the
``conductor'' of the ``orchestra'' of supply chain
members. Forecasting the future demand is, therefore, crucial for all
logistics activities. Special problems occur in cases like ours, when
the majority of inventory items is not replenished, because the
relationship between lead times and fashion cycles makes
replenishments simply impossible. The resulting ``textile pipeline''
has strong interdependencies between marketing, procurement, and
logistics.

The business model of our real world problem bases on a strict cost
leadership strategy with sourcing in low cost countries, either East
Asia or South East Europe. The transportation time is between one and
six weeks, economies of scale are achieved via large orders.

\subsection{Internal stock turnover of pre-packs}
\label{sec:real-world-problem:stock-turnover}

The material flow in our problem is determined by a central
procurement for around 1\,200 branches. All items are delivered
from the suppliers to a central distribution center, where a so called
``slow cross docking'' is used to distribute the items to the
branches. Some key figures may give an impression of the situation:
32\,000 square meters, 80 workers, 30\,000 tons of garment in 10
million lots per year.
%%, average stock level 8 million euro. 
Each branch
is delivered once a week with the help of a fixed routing system. This
leads to a sound compromise between inventory costs and costs of stock
turnover.

There are two extreme alternatives for the process of picking the
items. The retailer can either work with one basic lot and deliver
this lot or integer multiples of it for every article to the branches,
or he develops individual lots for the shops. At the beginning of our
analysis most of the articles were picked in one basic lot, but more
than 40 other constellations were used additionally. The costs for
picking the items are not relevant for the following analysis because
only minor changes of internal processes are necessary. Only in later
phases of our project a detailed analysis of the picking costs will be
needed.

\subsection{The problem under consideration}
\label{sec:real-world-problem:problem-description}

Recall that the stock turn-over is accelerated by ordering pre-packs of every
product, i.e., a package containing a specified number of items of each size. We
call the corresponding vector with a non-negative integral entry for every size
a \emph{lot-type}.

In this environment, we focus on the following problem: Given
historic delivery data (in terms of pre-packed lots of some lot-type
for each branch) and sales data for a group of products for each
branch, determine for each branch a new lot-type with the same number
of items that meets the relative demand for sizes more accurately.

In particular, find out from historic \emph{sales} data some
information about the relative \emph{demands} for sizes.  We stress
the fact that stockout substitution in the data can not be neglected
since replenishment does not take place (unsatisfied demand is lost
and does not produce any sales data).  We also stress the fact that we
are not trying to improve the number of items delivered to each branch
but the distribution of sizes for each branch only.

\section{Some real-data analysis evaluating seemingly obvious approaches}
\label{sec:data-analysis}

To deal with the problem described in
Section~\ref{sec:real-world-problem} our partner has provided us some
historic sales data for approximately $1\,200$~branches over a time
period of 12~months. The task was to forecast the size distribution of
the future demand for each branch. A possible size distribution is
depicted in Figure~\ref{tab:size_distribution}.

\begin{table}[htbp]
  \begin{center}
    \includegraphics[width=5cm]{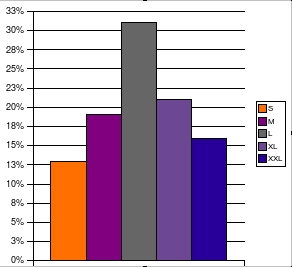} % aus Data/abweichung.xls
  \end{center}
  \caption{A possible size distribution of the demand}
  \label{tab:size_distribution}
\end{table}

The most obvious way to determine a size distribution as in
Figures~\ref{tab:size_distribution} is to count the number of sold items per
size and divide by the total number of sold items. Here we have some freedom to
choose the day of the sale where we measure the amounts. We have to balance two
competing facets. An early measurement may provide sales figures which are
statistically too small for a good estimate while a late measurement may suffer
from unsatisfied demand that is not present in the sales data. 

The business strategy of our partner implies to cut prices until all items are
sold. So, a very late measurement would only estimate the supply instead of the
demand. As there is no expert knowledge to decide which is the \textit{optimal}
day of sales to count the amounts and estimate the size distribution we have
applied a statistical test to measure the significance of the size distributions
obtained for each possible day of counting the sold items. 

Given a data set $D$, a day of sales $d$, and a size $s$ let $\phi_{s,d}(D)$ be
the estimated demand for size $s$ measured on day $d$ as described above. We
normalize so that $\sum\limits_s \phi_{s,d}(D)=1$ for each day and each data
set. Our statistical test partitions the original data set $D$ randomly into two
disjoint data sets $D_1$ and $D_2$. Naturally we would not trust a forecast
$\phi$ whenever $\phi_{s,d}(D_1)$ and $\phi_{s,d}(D_2)$ are too far apart.
Statistically speaking $\phi$ would not be robust in that case. To measure more
precisely how far apart $\phi_{s,d}(D_1)$ and $\phi_{s,d}(D_2)$ are, we define
the discrepancy $\delta_d$ as
$\delta_d(D_1,D_2):=\sum\limits_s\left|\phi_{s,d}(D_1)-\phi_{s,d}(D_2)\right|$.
In Figure \ref{tab:discrepancy-result} we have depicted the average discrepancy
over all branches for the first $60$ days of sale.

\begin{table}[htbp]
  \begin{center}
    \includegraphics[width=\linewidth]{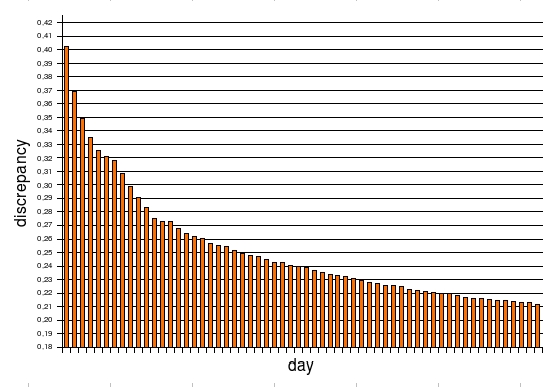} % aus Data/abweichung.xls
  \end{center}
  \caption{Average discrepancy $\delta_d$ over all branches}
  \label{tab:discrepancy-result}
\end{table}

For our practical application the significance of $\phi$ is too low
for all possible days of measurement. We remark that in our examined
data set the discrepancy $\delta_d$ tends to $0.19$ as $d$ tends to
the end of the sales period. An explanation why this obvious approach
does not work in our case is due to the small sales numbers and the
interference of product attractivity and price cutting strategies.

Even when we tried to use only use ordinal information generated from
an estimated size distribution (some size is too scarce, some size too
ample), we encountered different results from different samples (same
size was scarce in one sample and ample in the other).

In the same manner we have checked some common parametric models
(mostly based on an estimation of a constant sales rate for each
individual product-size pair) for demand forecasting on historic sales
data from literature. None of them gives significant information of
the size distribution of the demand of our data set, since the sales
rates vary drastically and depend more on the attractivity of the
products than on the sizes. For the details we refer to two diploma
theses \cite{Kouris:Diplom:2007,Wopperer:Diplom:2007} from our group.

\section{The Top-Dog-Index (TDI)}
\label{sec:TDI}

In the previous section we learned that in our application, first, we cannot
trust the common parametric models for the demand distribution and, second,
measuring a size profile directly from sales data may lead to more or less
random decisions. The main reason for the latter is the former and the
interference between attractiveness of offered products and compatibility of
offered sizes. Because of this, the stockout saturation of sales data happens
at almost all times during the sales period, and thus aggregating the sales data
of different products yields no reasonable information.

Our new idea throws overboard the desire to estimate an absolute size profile
of the demand in every branch.  Instead, we try to define a measure for the
scarcity of sizes during the sales process that can be estimated from historic
sales data in a stable way.

The following thought experiment is the motivation for our distribution free
measure: Consider a product, for which in a branch all sizes are sold out at the
very same day.  This can be regarded as the result of an ideal balance between
sizes in the supply.  Our measure tries to quantify the deviation of this
ideal situation in historic sales data.  How can this be done?  In the
following, we extract data of a new type from the sales process.

Fix a delivery period~$\Delta := [0, T]$ from some day in the past normed to
day~$0$ to day~$T$. Let $B$ be the set of all branches that are operating in
time interval~$\Delta$, and let $P$ be the set of all products in a group
delivered in time interval~$\Delta$ in sizes from a size set~$S$. We assume that
in each branch the product group can be expected to have homogeneous demand for
sizes throughout the time period. Fix $b \in B$. For each $p \in P$ and for each
$s \in S$ let $\theta_b(p, s)$ be the stockout-day of size~$s$ of product~$p$,
i.e., the day when the last item of $p$ in size~$s$ was sold out in branch~$b$.

Fix a size $s \in S$. Our idea is now to compare for how many products $p$
size~$s$ has the earliest stockout-day~$\theta(p, s)$ and for how many products
$p$ size~$s$ has latest stockout-day~$\theta(p, s)$. These numbers have the
following interpretation: If for many more products the stockout-day of the
given size was first among all sizes, then the size was scarce. If for many more
products the stockout-day of the given size was last among all sizes, then the
size was ample.

In order to quantify this, we use the following approach.  (In fact, it is not
too important how we exactly quantify our idea, since we will never use the
absolute quantities for decision making; we will only use the quantities
relative to each other.)

\begin{Definition}[Top-Dog-Index]
  Let $s \in S$ be a size and $b$ be a branch.  
  
  The \emph{Top-Dog-Count} $\TDC{b}{s}$ for $s$ in $b$ is defined as
  \begin{equation}
    \TDC{b}{s} := \bigl\lvert \{\; p \in P \;|\; \theta_b (p, s) = \min_{s' \in S}
    \theta_b (p, s')\; \} \bigr\rvert
  \end{equation}
  and the \emph{Flop-Dog-Count} $\FDC{b}{s}$ in $b$ is defined as
  \begin{equation}
    \FDC{b}{s} := \bigl\lvert \{\; p \in P \;|\; \theta_b (p, s) = \max_{s' \in S}
    \theta_b (p, s')\; \} \bigr\rvert.
  \end{equation}

  Moreover, for a fixed dampening parameter $C > 0$ let
  \begin{equation}
    \TDI{b}{s} := \frac{\TDC{b}{s} + C}{\FDC{b}{s} + C}
  \end{equation}
  be the \emph{Top-Dog-Index (TDI)} of Size~$s$ in~Branch~$b$.
\end{Definition}

In the data of this work, we used $C = 15$.

\subsection{Statistical significance of the Top-Dog-Index}
In a similar way as in Section \ref{sec:data-analysis} we want to
analyze the significance of the proposed Top-Dog-Index. Since this
method is supposed to be applied to a real business case, we analyze
the statistical significance in more detail. Instead of two data sets
$D_1$ and $D_2$ as in Section \ref{sec:data-analysis} we utilize seven
such sets $D_i$. Therefore, we assign a random number in $\{1,2,3,4\}$
to each different product. The sets are composed of the data of
products where the corresponding random number lies in a
characteristic subset of $\{1,2,3,4\}$. See Table \ref{tab:test_sets}
for the assignment. For the interpretation we remark that the pairs
$(D_1,D_2)$, $(D_3,D_4)$, and $(D_5,D_6)$ are complementary. The whole
data set is denoted by $D_7$.

\begin{table}[htbp]
  \begin{center}
    \begin{tabular}{rc}
      \hline
      $D_1$ & $\{1,2\}$\\
      $D_2$ & $\{3,4\}$\\
      $D_3$ & $\{1,3\}$\\
      $D_4$ & $\{2,4\}$\\
      $D_5$ & $\{3\}$\\
      $D_6$ & $\{1,2,4\}$\\
      $D_7$ & $\{1,2,3,4\}$\\
      \hline
    \end{tabular}
  \end{center}
  \caption{Assignment of test sets}
  \label{tab:test_sets}
\end{table}

Since the Top-Dog-Index is designed to provide mainly ordinal
information, we have to use another statistical test to make sure that
it yields some significant information. Let $\TDI{b}{s,D_i}$ denote
the Top-Dog-Index in Branch~$b$ of Size~$s$ computed from the data in
Data Set~$D_i$. We find it convincing to regard the ordinal
information generated by the Top-Dog-Index as robust whenever we have
\begin{equation}
  \TDI{b}{s,D_i} \gg \TDI{b}{s',D_i} \iff \TDI{b}{s,D_j} \gg \TDI{b}{s',D_j}
  \label{eq:tdi-order}
\end{equation}
for each pair of sizes $s, s'$ and each pair of data sets
$D_i,D_j$. In words: the order of Top-Dog-Indices of various sizes
does not change significantly when computed from a different sample.
The following is a sufficient condition for this to happen:
\begin{equation}
  \frac{\TDI{b}{s,D_i}}{\sum_j\TDI{b}{s,D_j}}
  \approx
  \mathit{const}_i
  \label{eq:tdi}
\end{equation}

Our first aim is to provide evidence that the $\TDI{b}{s}$ values are
robust measurements in this sense.  There is a nice way to look at
equation \ref{eq:tdi} graphically. Let us fix a size $s$. For each
branch $b$ let us plot a column of the relative values
$\frac{\TDI{b}{s,D_i}}{\sum_j\TDI{b}{s,D_j}}$ for all branch-size
combinations and for all~$i$.  The columns corresponding to the same
branch-size combination but different samples are stacked on top of
each other.  This way, each colum stack has height one in total, and
the size relations of the measurements based on the different samples
can be assessed right away.  To show the robustness of the mean and
the median, we have added an additional column for each of them: The
median of the seven estimates corresponds to the height of the
top-most column, the mean to the height of the second-to-top-most
column.  The goal now is not to read off the exact values, but to get
an intuitive impression how heavy an estimate depends on the data
subset it was computed from.

\begin{figure}[htbp]
  \begin{center}
    \includegraphics{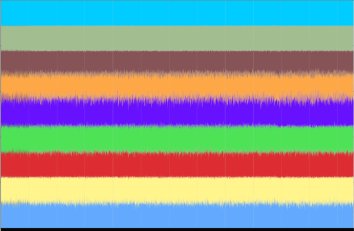}
  \end{center}
  \caption{Relative distribution of the Top-Dog-Index on seven data
    subsets for all branch-size combinations; the two top-most columns
    are median and mean, resp.}
  \label{tab:tdi-significance}
\end{figure}

The Top-Dog-Indices are plotted this way in
Figure~\ref{tab:tdi-significance}.

\begin{figure}[htbp]
  \begin{center}
    \includegraphics{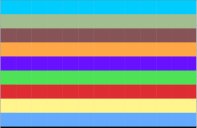}\quad
    \includegraphics{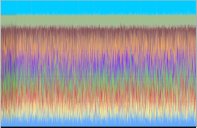} % aus NKD-Bericht
  \end{center}
  \caption{Relative distribution of deterministic and random numbers
    on seven data subsets and branches; the two top-most columns are
    median and mean, resp.}
  \label{tab:deterministic_vs_random}
\end{figure}

In order to provide some intuitively clear reference data to compare
to the plot of Figure~\ref{tab:tdi-significance}, we present the
corresponding plots for the two extreme cases of deterministic numbers
(i.e., $\TDI{b}{s,D_i}=\TDI{b}{s,D_j}$ for all $i$ and $j$) and
totally random numbers in Figure~\ref{tab:deterministic_vs_random}. In
the complete deterministic case the areas of same color form perfect
rectangles. In the random case the areas of same color corresponding
to the data subsets form zig-zag lines; the median has fewer zig-zag
than the mean, but both are quite stable because the random numbers
are all from the same distribution. 

It is immediately obvious that the plot of
Figure~\ref{tab:tdi-significance} looks more like the plot in the
deterministic case as the plot in the random case.  It is interesting
to note that the dampening parameter in the computation of the TDI
does indeed influence the amount of noise in the plot but had almost
no influence on the order of TDI values, which is what we indend to
use.

\begin{figure}[htbp]
  \begin{center}
    \includegraphics{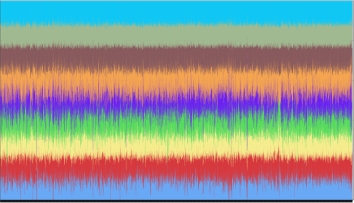} % aus NKD-Bericht
  \end{center}
  \caption{Relative distribution of sales up to Day~$0$ on on seven
    data subsets for all branch-size combinations; the two top-most
    columns are median and mean, resp.}
  \label{tab:estimate-significance-day-0}
\end{figure}

\begin{figure}[htbp]
  \begin{center}
    \includegraphics{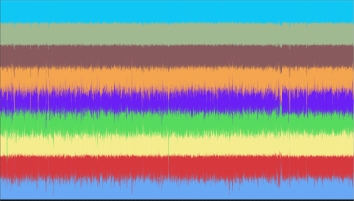}
  \end{center}
  \caption{Relative distribution of sales up to Day~$12$ on on seven
    data subsets for all branch-size combinations; the two top-most
    columns are median and mean, resp.}
  \label{tab:estimate-significance-day-12}
\end{figure}

Recall that we have analyzed the direct estimation of size
distributions for each branch by looking at two complementary samples
of the data set.  In order to be able to directly compare the results
of that attempt to Figure~\ref{tab:tdi-significance}, we plot the size
distributions for the same seven data subsets.  This is shown in
Figure~\ref{tab:estimate-significance-day-0} for an estimate from the
sales up to Day~$0$ (the first day in the sales period) and in
Figure~\ref{tab:estimate-significance-day-12} for an estimate from the
sales up to Day~$12$.  

It can be seen that Day~$0$-estimates are extremely dependent on the
data subset even if only the relative values are taken into account,
i.e., the estimates are not robust even in the weak sense measured in
the plot.  The Day-$12$ estimates are more robust, but not even close
to what the TDI achieves in Figure~\ref{tab:tdi-significance}.
Moreover, as we said, those estimates are already quite close to the
supply because of unsatisfied demand, and will therefore fail to
measure the size distribution of the demand.

%% We remark that Figure~\ref{tab:tdi-significance} corresponds to a
%% single size; yet, the plots for the other sizes or further data sets
%% are rather similar.

Generating a robust statistics for a set of data is, of course, easy:
just assign the same deterministic number to each object.  This is not
what we want since this does not carry any useful information.  We
claim that the Top-Dog-Index exhibits differences between demands for
sizes in each branch individually.  That this is indeed the case,
follows from the Top-Dog-Indices that we encountered in the field
study documented below.  But let us first discuss which actions we
could take to improve the size distributions of the branches' supply.

So far we have argued that the Top-Dog-Index produces size-related
information in a robust way, while other methods fail (at least for
our given real-world data). In the next section we describe a
procedure to harmonize demand and supply with respect to the size
distribution. In Section~\ref{sec:blind-study} we provide evidence via
a real-world blind study that this procedure helps to raise the gross
yield in reality; thus, the Top-Dog-Index is correlated to the size
distribution of the demand.

\section{The heuristic size optimization procedure based on the TDI}
\label{sec:heuristic-optimization}

Interesting for us is not the absolute TDI of a size but the TDIs of
all sizes in a branch compared to each other, i.e., the ordinal
information implied by the TDIs.  The size with the maximal TDI among
all sizes can be interpreted as the scarcest (the one that was sold
the fastest) size in that branch; the size with the minimal TDI among
all sizes can be interpreted as the amplest size in that branch.
Of course, we have the problem of deciding whether or not a maximal
TDI is significantly larger than the others.  Since the absolute
values of the TDIs have no real meaning we did not even try to assess
this issue in a statistically profound way.

Our point of view is again that absolute forecasting is too much to be
asked for.  Therefore, we resort to a dynamic heuristic optimization
procedure: sizes with ``significantly'' large TDIs (Top-Dog-Sizes)
should receive larger volumes in future deliveries until their TDIs
do not improve anymore, while sizes with ``significantly'' small TDIs
(Flop-Dog-Sizes) need smaller supplies in the future.  Whenever this
leads to oversteering, the next TDI analysis will show this, and we go
back one step.  This is based on the assumption that the demand for sizes
does not change too quickly over time.  If it does then optimization
methods based on historic sales data are useless anyway.

Let us describe our size distribution optimization approach in more detail.

We divide time into delivery periods (e.g., one quarter of a year).  We assume
that the sales period of any product in a delivery period ends at the end of
the next delivery period (e.g., half a year after the beginning of the
delivery period).  Recall that a size distribution in the supply of a branch
is given by a pre-pack configuration: a package that contains for each size a
certain number of pieces of a product (compare
Section~\ref{sec:real-world-problem:stock-turnover}).

We want to base our delivery decisions for an up-coming period on 
\begin{itemize}
\item the pre-pack configuration of the previous period and
\item the TDI information of the previous period giving us the deviation from
  the ideal balance
\end{itemize}
According to given restrictions from the distribution system, we
assume that only one pre-pack configuration per branch is allowed.  We
may use distinct pre-pack configurations for different branches,
though.

Since we are only dealing with the size distribution of the total
supply but not with the total supply for a branch itself, the total
number of pieces in a pre-pack has to stay constant.  Since the TDI
information only yields aggregated information over all products in
the product group, all products of this group will receive identical
pre-pack configurations in the next period, as desired.

In order to adjust the supply to the demand without changing the total number
of pieces in a pre-pack, we will remove one piece of a Flop-Dog-Size
from the pre-pack and add one piece of a Top-Dog-Size instead.
At the end of the sales period (i.e., at the end of the next delivery period),
we can do the TDI-analysis again and adjust accordingly.

Given the usual lead times of three months this leads to a heuristics
that reacts to changes in the demand for sizes with a time lag of nine
month to one year.  Not exactly prompt, but we assume the demand for
sizes to be more or less constant over longer periods of time.

The most interesting question for us was how much, in practice, could
be gained by performing only one step of the heuristics explained
above.

\section{A real-world blind study}
\label{sec:blind-study}

In this section we describe the set-up and the results of the blind study
carried out by our partner.  A summary of parameters can be found in
Table~\ref{tab:summary}.

\begin{table}[htbp]
  \begin{center}
    \begin{tabular}{ll}
      \hline
      Test period :           & April through June 2006 (3~Months)\\
      Data collection period: & April through September 2006 (6~Months)\\
      Branches:               & $20$~branches with unbalanced TDIs\\
      & randomly classified into $10$~test and $10$~control branches\\
      Pieces of merchandise:  & approx.~$4\,000$ pieces for all test and control branches\\
      Merchandise value:      & approx.~$30\,000$\,\euro\\
      \hline
    \end{tabular}
  \end{center}
  \caption{Summary of parameters of the blind study}
  \label{tab:summary}
\end{table}

\subsection{Selection of branches}
\label{sec:blind-study:selection-of-branches}
  
A reasonable selection of branches for a test and a control group had
to meet essentially three requiremens: first, only those branches
should be chosen whose TDI indicated that, in the past, the supply by
sizes did not meet the demand for sizes; secondly, no branch should be
chosen, where other tests were running during the test period;
thirdly, the assignment of branches to test- and control group should
be completely random.  The reason for the third aspect was that this
way all other influences on the gross yield than the selection of
sizes would appear similarly in both the test group and the control
group and, thus, would average out evenly.

We suggested a set of 50~branches with interesting Top-Dog-Indices to
our partner. Out of these 50~branches, our partner chose 20 branches
where a potential re-packing of pre-packs would be possible.  This set
of 20~branches was fixed as the set of branches included in our blind
study.

After that, a random number between 0 and~1 was assigned to each of
the 20~branches.  The 10~branches with the smallest random numbers
were chosen to be the test group, the rest was taken as the control
group.

\subsection{Handling of pre-packs}
\label{sec:blind-study:pre-pack-handling}

Next, we had to specify the modifications to the size distribution in
pre-packs on the basis of the TDI information.  It turned out that
additional side constraints had to be satisfied: Whenever a product
would appear in an advertisement flyer, the pre-pack had contained at
least one piece in each of the four main sizes S, M, L, and~XL.  That
means, although sometimes the TDI suggested that S was the amplest
size, we could not remove the only piece in~S from the pre-pack.  We
removed a piece of the second amplest size (M or L, of which there
were two in the unmodified pre-packs) instead.  To all branch
deliveries, one additional piece of~XL was packed, because this was
the scarcest size in every branch in the test group.  This way, the
total number of pieces was unchanged in every pre-pack, as suggested
in Section~\ref{sec:heuristic-optimization}.

Since all orders had been placed well before the decision to conduct a
blind study, our partner re-packed all pre-packs for the test branches
according to Table~\ref{tab:pre-pack-modifications}.

\begin{table}
  \begin{center}
    \begin{tabular}{c|cc|cc}
      & \multicolumn{2}{c|}{\textbf{special}}
      & \multicolumn{2}{c}{\textbf{no}} \\
      \textbf{test group}
      & \multicolumn{2}{c|}{\textbf{advertising}} 
      & \multicolumn{2}{c}{\textbf{advertising}} \\
      Branch & remove & add & remove & add \\
      \hline
       1 & L & XL & L & XL\\
       2 & M & XL & M & XL\\
       3 & L & XL & S & XL\\
       4 & L & XL & S & XL\\
       5 & M & XL & M & XL\\
       6 & M & XL & M & XL\\
       7 & M & XL & M & XL\\
       8 & L & XL & S & XL\\
       9 & M & XL & M & XL\\
      10 & M & XL & M & XL
    \end{tabular}
  \end{center}
  \caption{How the pre-packs were modified in the test group}
  \label{tab:pre-pack-modifications}
\end{table}

\subsection{Time frame}
\label{sec:blind-study:time-frame}

The test included two relevant time periods: the first period from
which the TDI data was extracted and the test period in which the
recommendations based on the TDI data were implemented for the test
group.

The TDI data was drawn from a delivery period of nine month (January
through September 2005) and a sales period of twelve month (January
through December 2005).

The test data was drawn from a delivery period of three months (April
through June 2006) and a sales period of six months (April through
September 2006).

\subsection{Data collection}
\label{sec:blind-study:data-collection}

In order to sort out contaminated data easily, our partner agreed to take
stock to check inventory data for correctness every month. To receive a good 
estimate for the financial benefit of the supply modification proposed by the 
procedure described in the previous section we had defined some criteria how 
to detect contaminated data automatically via a computer program.

\subsection{Data analysis}
\label{sec:blind-study:data-analysis}

%%%%%%%%%%%%%%%%%%%%%%%%%%%%%%%%%%%%%%%%%%%%%%%%%%%%%%%%%%%%%%%%%%
%% Mean TDIs before:
%%%%%%%%%%%%%%%%%%%%%%%%%%%%%%%%%%%%%%%%%%%%%%%%%%%%%%%%%%%%%%%%%%

\begin{figure}[tbp]
  \begin{center}
    \includegraphics[width=0.85\linewidth]{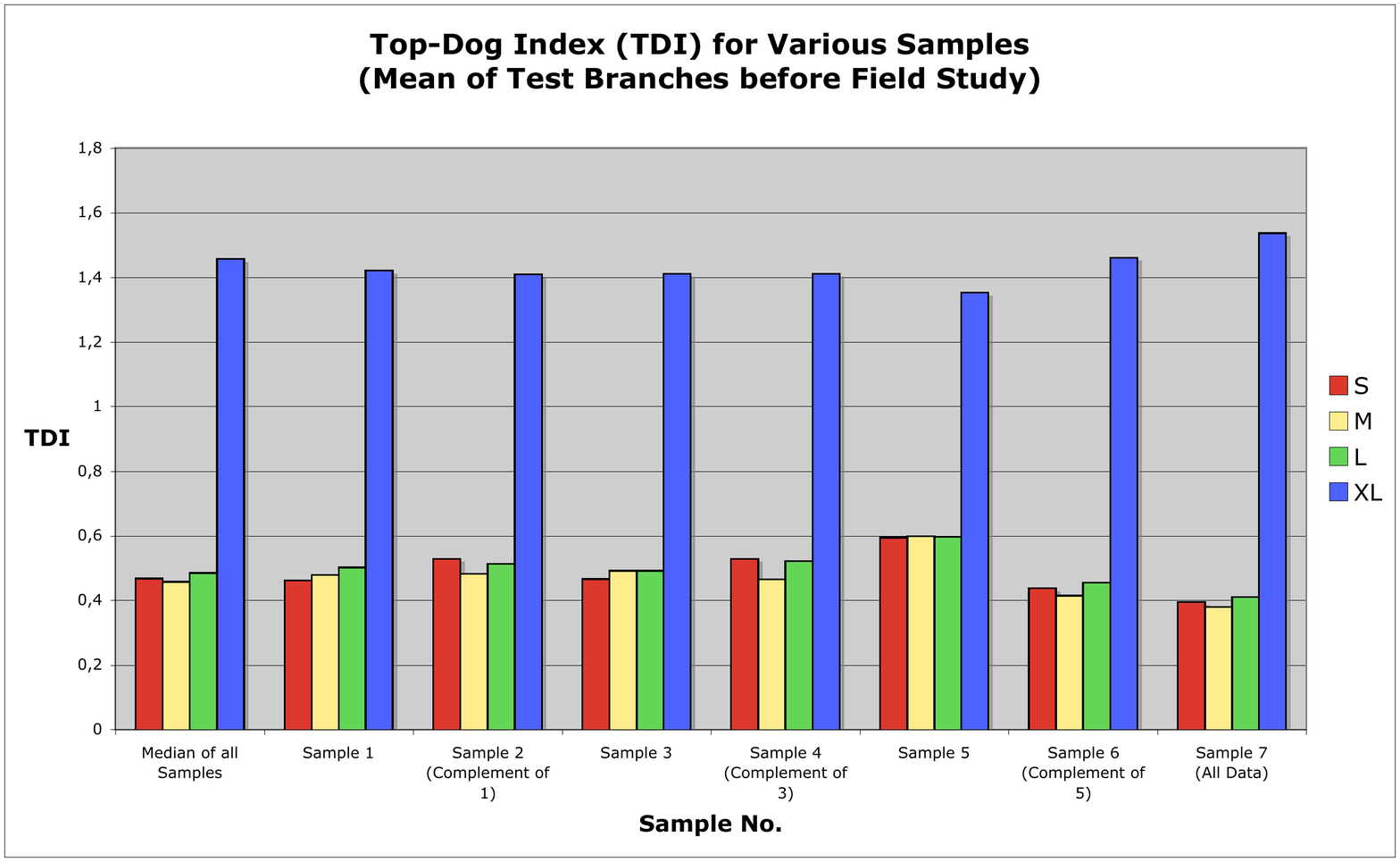}
  \end{center}
  \caption{TDI in the test branches from historic sales data}
  \label{tab:TDI-before-test}
\end{figure} 

\begin{figure}[tbp]
  \begin{center}
    \includegraphics[width=0.85\linewidth]{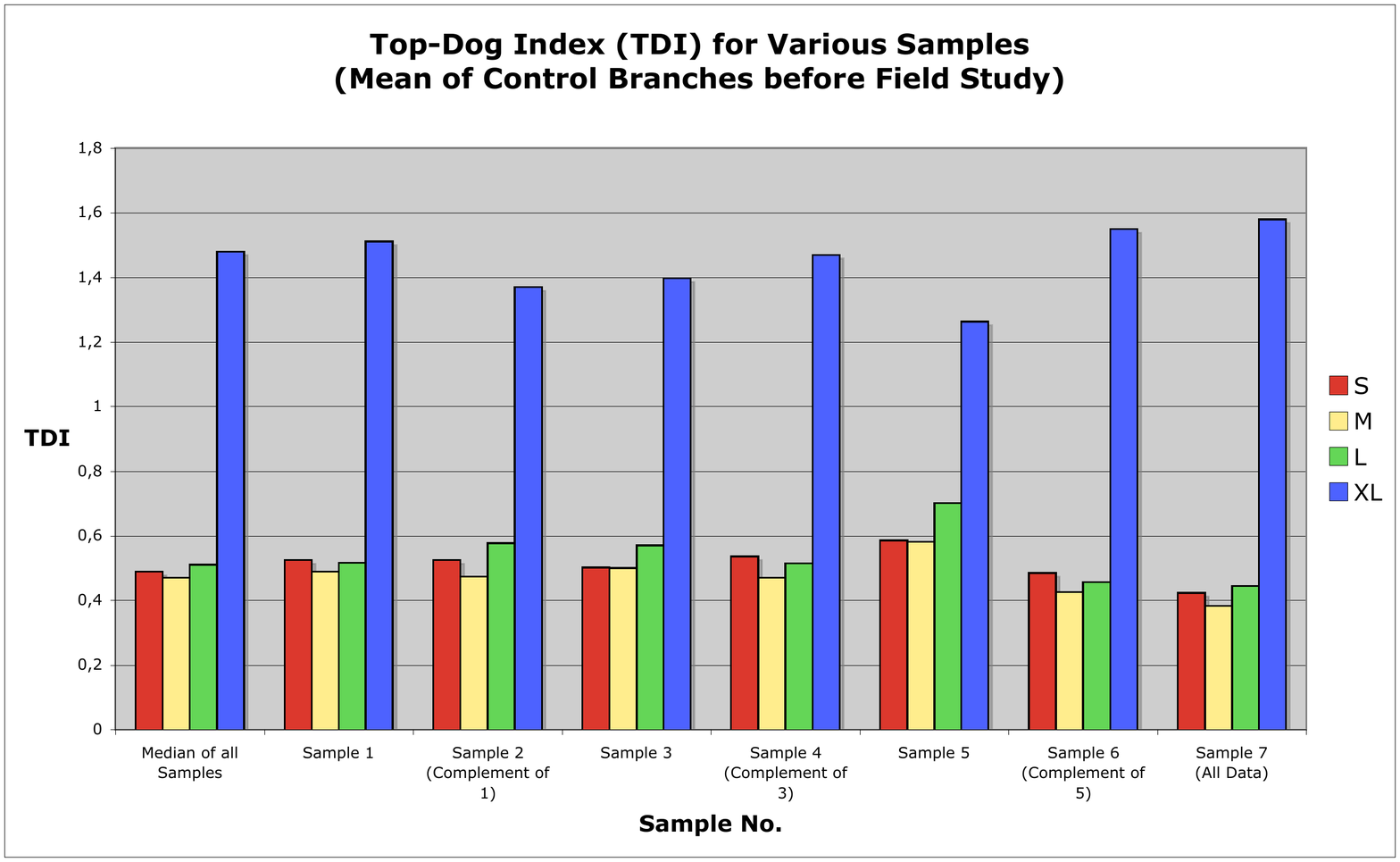}
  \end{center}
  \caption{TDI in the control branches from historic sales data}
  \label{tab:TDI-before-control}
\end{figure}

%%%%%%%%%%%%%%%%%%%%%%%%%%%%%%%%%%%%%%%%%%%%%%%%%%%%%%%%%%%%%%%%%%
%% TDIs before of Individual Branches:
%%%%%%%%%%%%%%%%%%%%%%%%%%%%%%%%%%%%%%%%%%%%%%%%%%%%%%%%%%%%%%%%%%

\begin{figure}[tbp]
  \begin{center}
    \includegraphics{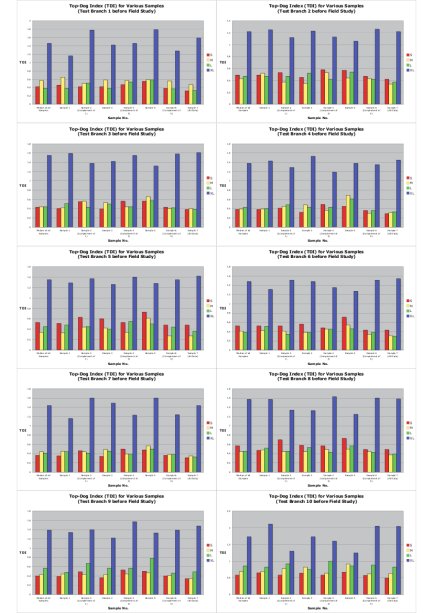}
  \end{center}
  \caption{Individual TDIs in test branches from historic sales data}
  \label{tab:TDI-before-test-all}
\end{figure}

\begin{figure}[tbp]
  \begin{center}
    \includegraphics{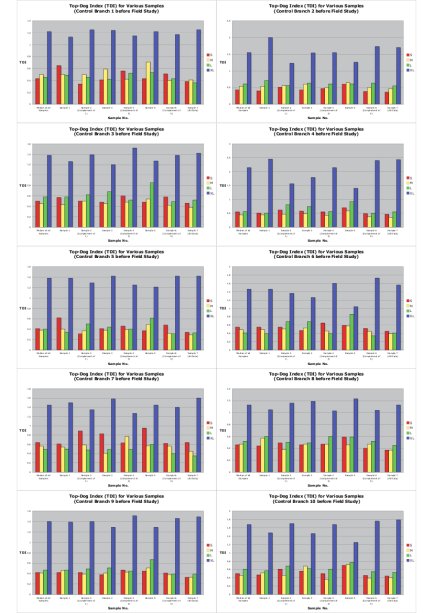}
  \end{center}
  \caption{Individual TDIs in control branches from historic sales data}
  \label{tab:TDI-before-control-all}
\end{figure}

The initial situation was given by the data set described in Section
\ref{sec:data-analysis}.  In Figure~\ref{tab:TDI-before-test} we have
depicted the initial Top-Dog-Indices for the test branches and in
Figure~\ref{tab:TDI-before-control} the initial Top-Dog-Indices for
the control branches.  Moreover, Figures \ref{tab:TDI-before-test-all}
and~\ref{tab:TDI-before-control-all} finally show the situations for
the individual branches.

The analysis of our field study was intended to answer the following
two main questions: are the Top-Dog-Indices better distributed in the
test branches than in the control branches, and, if yes, does this
have a significant monetary impact?  

To investigate the latter, we had to analyze some monetary
variables. The most important monetary indices for our partner are the
gross yield and the last price. The gross yield directly shows how
much turnover was lost using a price cutting strategy to sell out all
items. The last price tells us how far one was forced to dump items
provoked by an inadequate size distribution of the supply. Since the
values of different products vary widely we only consider relative
values. So, the gross yield is given by the ratio of realized turnover
and theoretic revenue without cutting the prices.

Since we had to deal with a large amount of lost or inconsistent data,
we have applied two ways of evaluating gross yield and last
price. Imagine that your data says that you have sold 10~items of a
given article in a given branch but that the supply was only
8~items. Or the other way round that the initial supply was 10~items,
during the time 8~items were sold, but all items are gone.  Both of
the described situations occurred significantly often in our data set.

Our first strategy to evaluate the given data was to ignore
inconsistent data. In the first case, 8~sale transactions are
consistent. For the remaining two items the corresponding supply
transaction is missing. So, we simply ignore these transactions. In
the other case we would ignore the supply of two items.

The alternative to ignoring inconsistent data is to estimate or
recover it from the rest of the data set. As an example, we would
simply assume that there was a supply of 10~items instead of 8~items
at the same price level in the first case. In the second case we would
assume that the remaining 2~items were also sold. Maybe they were
shoplifted, some sort of selling for a very cheap prize. So, we need an
estimation for the selling prize of the two missing items. Here, we
have used the last selling price over all sizes for this product in
this branch as an estimate. 

Neither evaluation method reflects reality exactly. Our hope was that
both estimations encompass the true values. At least our partner accepts
both values as a good approximation of reality. The truth may be
somewhere in between both values.

\subsection{Results}
\label{sec:blind-study:results}

\begin{figure}[tbp]
  \begin{center}
    \includegraphics[width=0.85\linewidth]{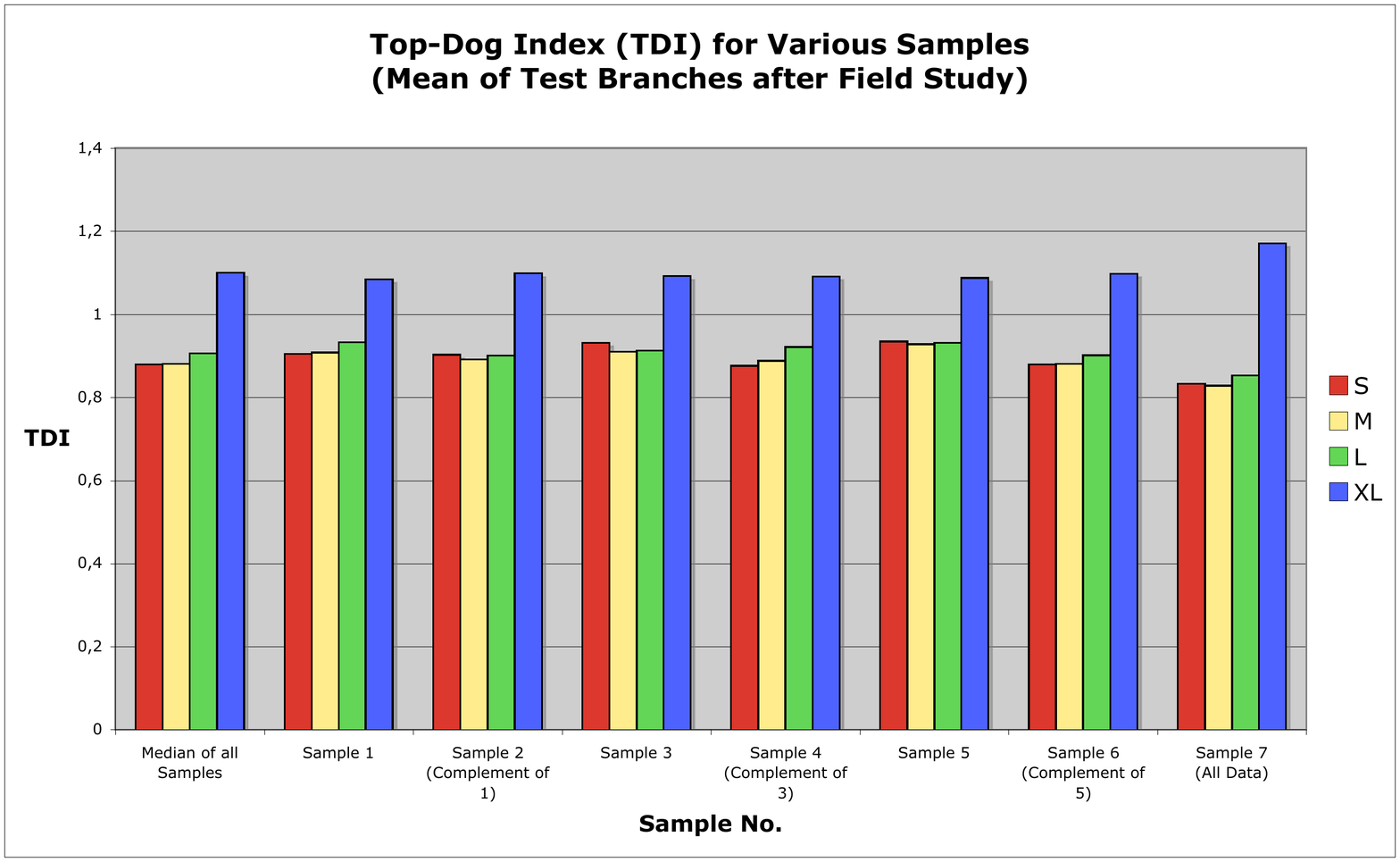}
  \end{center}
  \caption{TDI in the test branches in test period}
  \label{tab:TDI-after-test}
\end{figure}

\begin{figure}[tbp]
  \begin{center}
    \includegraphics[width=0.85\linewidth]{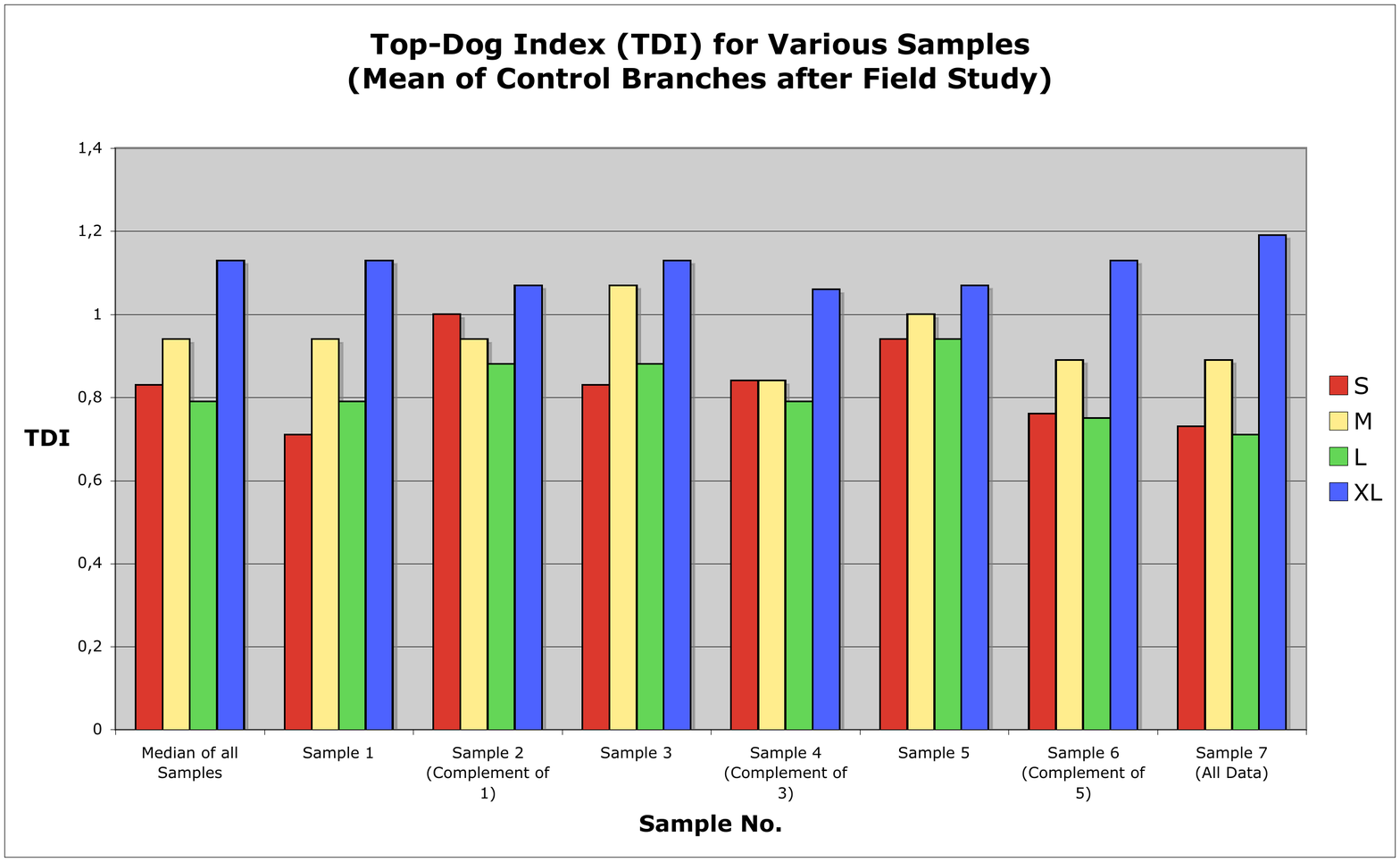}
  \end{center}
  \caption{TDI in the control branches in test period}
  \label{tab:TDI-after-control}
\end{figure}

\begin{figure}[tbp]
  \begin{center}
    \includegraphics{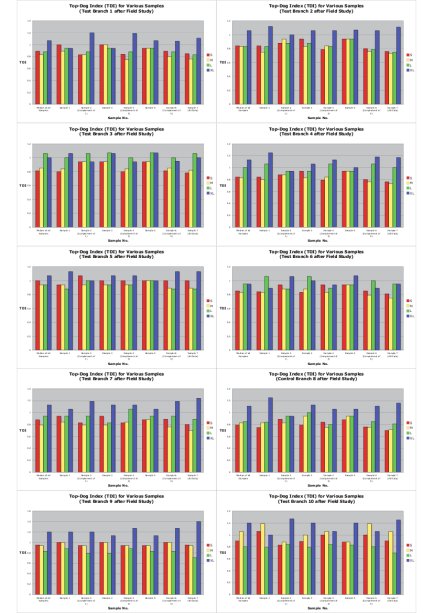}
  \end{center}
  \caption{Individual TDIs in test branches in test period}
  \label{tab:TDI-after-test-all}
\end{figure}

\begin{figure}[tbp]
  \begin{center}
    \includegraphics{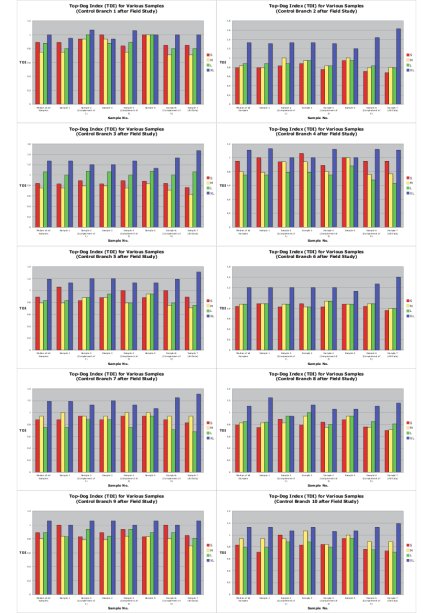}
  \end{center}
  \caption{Individual TDIs in control branches in test period}
  \label{tab:TDI-after-control-all}
\end{figure}

We have depicted the new Top-Dog-Indices after applying our proposed
repacking in Figure \ref{tab:TDI-after-test} for the test branches and
in Figure \ref{tab:TDI-after-control} for the control branches.
Figures \ref{tab:TDI-after-test-all}
and~\ref{tab:TDI-after-control-all} finally show the results for the
individual branches.

%%%%%%%%%%%%%%%%%%%%%%%%%%%%%%%%%%%%%%%%%%%%%%%%%%%%%%%%%%%%%%%%%%
%% Mean TDIs after:
%%%%%%%%%%%%%%%%%%%%%%%%%%%%%%%%%%%%%%%%%%%%%%%%%%%%%%%%%%%%%%%%%%

%%%%%%%%%%%%%%%%%%%%%%%%%%%%%%%%%%%%%%%%%%%%%%%%%%%%%%%%%%%%%%%%%%
%% TDIs after of Individual Branches:
%%%%%%%%%%%%%%%%%%%%%%%%%%%%%%%%%%%%%%%%%%%%%%%%%%%%%%%%%%%%%%%%%%

We can see in Figure~\ref{tab:TDI-before-control-all}
and~\ref{tab:TDI-after-control-all} that the absolute Top-Dog-Indices
are not as constant over seasons as we would wish.  However, the
induced order on sizes seems quite stable accross seasons, and this is
all we try to exploit.

We can see that it is rather hard to compare the Top-Dog-Indices of
the same branches before and after the blind study. The situation on
the real market almost never stays constant over time. There are so
many influences not considered in our study that it would have been a
bad idea to measure a possible raise of earnings directly.  For this
reason, the simultaneous observation of a test group and a control
group makes all outer effects appear in both. 

Comparing Figure~\ref{tab:TDI-after-test} and
Figure~\ref{tab:TDI-after-control} based on the same time period, it
appears that the Top-Dog-Indices of the test group have improved more.

More specifically, in Figures~\ref{tab:TDI-after-test-all}
and~\ref{tab:TDI-after-control-all} we see that in some of the test
branches (especially, Test Branches 5 and~6), Size~XL is no longer too
scarce, while it remains too scarce in other branches.  Moreover, on
average over all test branches, the Top-Dog-Indices of Sizes S, M,
and~L are better balanced, while in the control branches the
corresponding Top-Dog-Indices differ.

While the former achievement might have been equally possible on the
basis of a statistics aggregated over all branches, it seems that the
latter result was made possible only by the branch dependent
information from the Top-Dog-Index, since different sizes were removed
from the pre-packs in favor of~XL.  Moreover, some of the test
branches still need fewer pieces in Size~XL, some test branches don't.
That is, in the next optimization step, the branch dependent
information becomes vital also for the supply with Size~XL.

But is an improved Top-Dog-Index really an improvement for the
business?  To answer this question, we have quantified the gross
yields and the last prices in the test group and the control group,
resp.

\begin{figure}[htbp]
  \begin{center}
    \includegraphics[width=7cm]{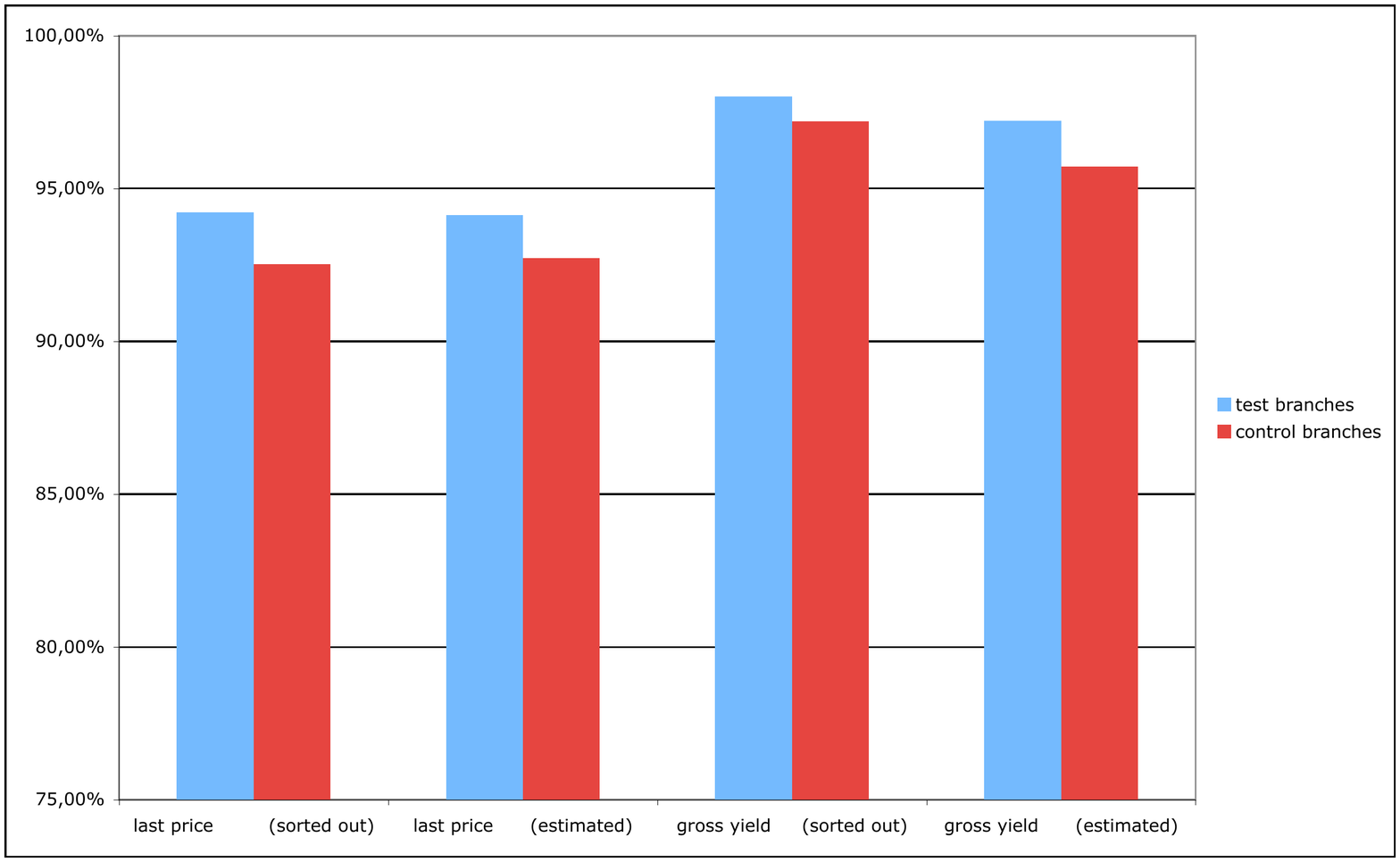}
  \end{center}
  \caption{Average last prices/gross yields per merchandise value}
  \label{tab:average-result}
\end{figure}

\begin{figure}[htbp]
  \begin{center}
    \includegraphics[width=7cm]{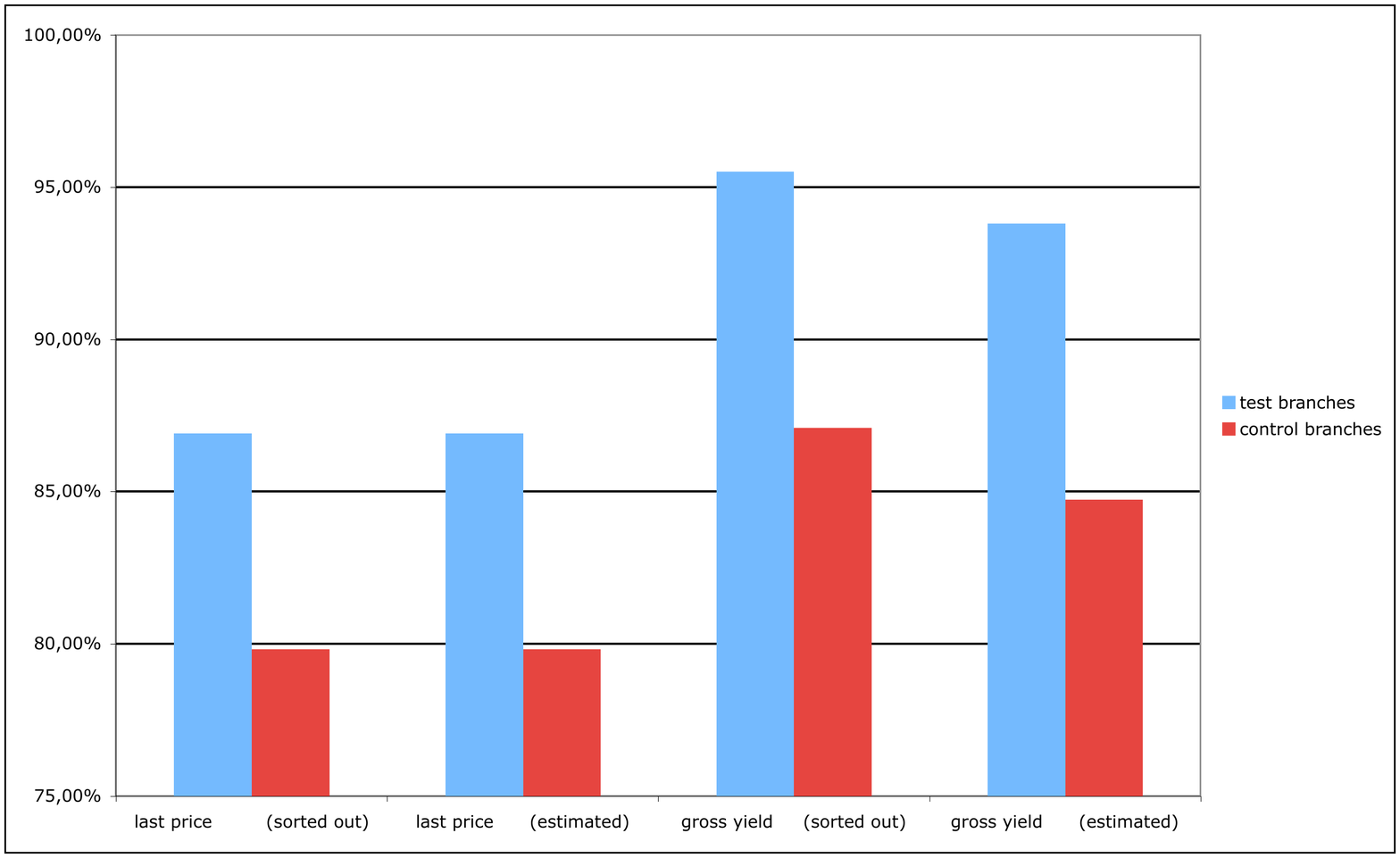}
  \end{center}
  \caption{Minimal last prices/gross yields per merchandise value}
  \label{tab:minimial-result}
\end{figure}

\begin{figure}[!h]
  \begin{center}
    \includegraphics[width=7cm]{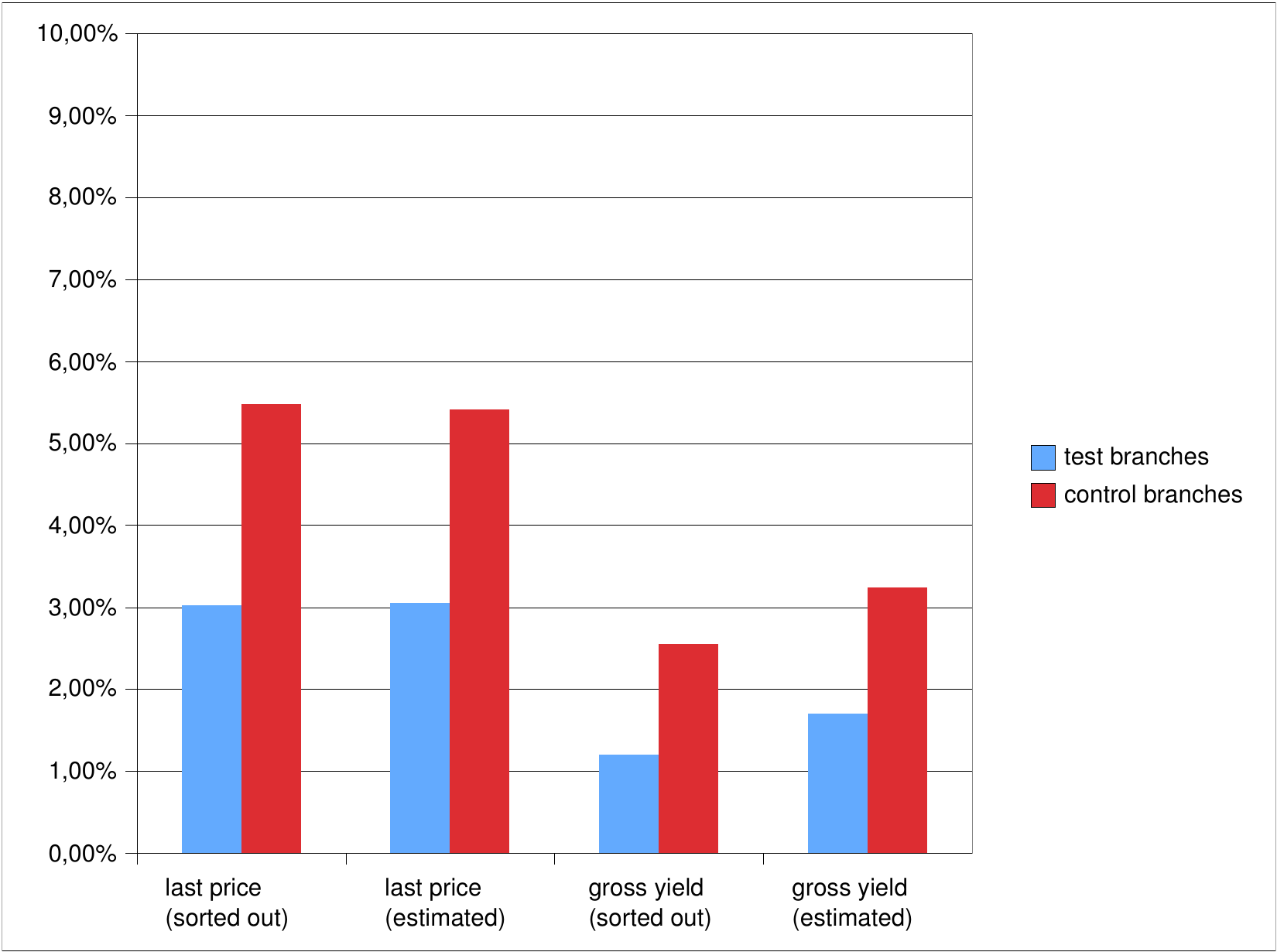}
  \end{center}
  \caption{Standard deviation of last prices/gross yields per merchandise value}
  \label{tab:deviation-result}
\end{figure}

In Figure~\ref{tab:average-result} we have compared the average values
of the gross yield and the last price for the control and the test
branches for both evaluation methods. The gross yield of the test
branches is 98.0\,\% using the ignore and 97.2\,\% using the estimate
method. For the control branches we have gross yields of 97.2\,\%
(ignore) and 95.7\,\% (estimate). This corresponds to improvements
of 0.85 and 1.5 percentage points, resp.

The improvement for the last price is even larger. The test branches
show a last price level of 94.2\,\% (ignore) and 94.1\,\%; the control
branches exhibit a last price level of 92.5\,\% (ignore) and 92.7\,\%
(estimate).  This corresponds to improvements of 1.7 and 1.4 percentage
points, resp.

The drastically improved results for the respective ``loser branches''
(see Figure~\ref{tab:minimial-result}) and the reduced standard
deviation (see Figure~\ref{tab:deviation-result}) in the data of the
test branches provides evidence for the fact that our procedure was
able to reduce the risk of a very low last price or a very low gross
yield in an individual branch.  This effect is desirable beyond the
better earnings, as very low last prices undermine the image of the
retailer.

\subsection{Statistical evidence of the improvement of the gross yield}

As in the previous sections we analyze our results concerning the
gross yield from the statistical point of view. Faced with the fact of
a widely varying gross yield over the branches with no appropriate
theoretic sales model, we have to restrict ourselves to
distribution-free statistics.

Therefore, we adapt the Wilcoxon rank sum test to our situation.  This
is a test method to find out whether or not two data sets are drawn
from the same distribution. We sort in increasing order the gross
yields of the $20$ branches participating in our blind study. We
associate the largest value with Rank~1, the second largest with
Rank~2, and so on.  Then we form the rank sums of the test branches
and the control branches, resp.  The more the rank sums differ, the
less likely is the event that our method did not influence the gross
yields/last prices at all.

For the Wilcoxon rank sum test, it is vital that we have partitioned
the 20~branches for the blind study independently at random into test
and control branches.

It is intuitively clear that a smaller rank sum for the gross
yield/last price is more likely if the corresponding expected values
are better. A lower rank sum can indeed occur by pure coincidence, but
the probability decreases with the rank sum. As an example, the rank
sum for the test branches regarding the gross yield measured by the
ignore method is~$89$. If we had not changed anything, the chance to
receive a rank sum of $89$ or lower would have been 12.4\,\%. So, we
have a certainty of 87.6\,\% that our proposed re-packing improved the
situation.  (More formally: the probability that the gross yields of
the test branches and the gross yields of the control branches stem
from the same distribution, i.e., nothing has changed systematically,
is at most 12.4\,\%.)

Now we consider different scenarios. Let $y_i(b)$ be the gross yield
measured with the method ignore of branch $b$ and $y_e(b)$ the gross
yield measured with the method estimate. By $i_c$ we denote the
scenario where we consider the values $y_i(b)$ for control branches
and the values $y_i(b)+\frac{c}{100}$ for test branches. Similarly, we
define the scenarios for $e_c$ utilizing $y_e(b)$ instead of
$y_i(b)$. In Table \ref{tab:rank_sums} we have given the rank sums and
the certainties of some scenarios.

\begin{table}[htbp] % Daten aus Data/Umsaetze_Feldversuch.xls
  \begin{center} % Wahrscheinlichkeiten mit Data/score.cpp berechnet (Zahl direkt neben Rangsumme zaehlt)
    \begin{tabular}{|r|rrrrrrrrrr|}
      \hline
      & $\!i_{-0.00}\!$ & $\!e_{-0.00}\!$ & $\!i_{-0.25}\!$ & $\!e_{-0.25}\!$ & $\!i_{-0.50}\!$ & $\!e_{-0.50}\!$ &
      $\!i_{-0.75}\!$ & $\!e_{-0.75}\!$ & $\!e_{-1.00}\!$ & $\!e_{-1.50}\!$\\
      \hline
      control group &  121 &  131 &  118 &  128 &  105 &  124 &   94 &  120 &  112 &  107 \\
      test group    &   89 &   79 &   92 &   82 &  105 &   86 &  116 &   90 &   98 &  103 \\
      certainty (\%)& 87.6 & 97.4 & 82.4 & 95.5 & 48.5 & 91.7 & 19.7 & 86.0 & 68.5 & 54.4 \\
      \hline
    \end{tabular}
  \end{center}
  \caption{Ranks sums and certainties of improvements of the gross yield}
  \label{tab:rank_sums}
\end{table}

How can we interprete these numbers? The first two columns of Table
\ref{tab:rank_sums} show that with a certainty of~87.6\,\% (ignore)
and 97.4\,\% (estimate) that our proposed modification increased the
expected gross yield. In Scenario~$i_{-0.25}$ we artificially decrease
the gross yield (ignore) values by 0.25 percentage points.  The
monetary value associated with this specific decrease can be
interpreted, e.g., as the implementation and consultancy costs of the
modification. So, by a look at Table \ref{tab:rank_sums} we can say
that with a certainty of 82.4\,\% our proposed modification yields an
improvement of the gross yield (ignore) by at least 0.25 percentage
points. 

%% Note that complete uncertainty would correspond to a certainty of
%% 50\,\%. We can see that certainty \textit{costs} something. In the
%% previous subsection we said that we have improved the gross yield
%% (ignore) by 0.5 percentage points in average. Table
%% \ref{tab:rank_sums} gives us a more refined analysis of the
%% improvement. But nevertheless the average improvement is a good
%% estimate for the real improvement since the rank sum test is a very
%% weak test that uses almost no information. The interested reader is
%% invited to check for himself how the rank sum test performs on a
%% deterministic improvement on some similar data.

%% \section{Implementation issues in the holistic view}
%% \label{sec:implementation-issues}

%% \subsection{Implications for the supply chain}
%% \label{sec:implementation-issues:implementation-supply-chain}

%% \subsection{Implications for the internal stock turnover}
%% \label{sec:implementation-issues:implementation-stock-turnover}

\section{Conclusion and outlook}
\label{sec:conclusion-and-outlook}

The distribution of fashion goods to the branches of a fashion
discounter must meet the demand for sizes as accurately as possible.
However, in our business case, an estimation of the relative demand
for apparel sizes from historic sales data was not possible in a
straight-forward way.  

Our proposal is to use the \emph{Top-Dog-Index (TDI)}, a measure that
yields basically ordinal information about what were the scarcest and
the amplest sizes in a product group in a historic sales period.  This
information was utilized to change the size distributions for future
deliveries by replacing one piece of the amplest size by a piece of
the scarcest size in every pre-pack (this can be seen as a subgradient
improvement step in an iterative size distribution heuristics based on
the TDI analysis).

Empirical evidence from a blind study with twenty branches (ten of
them, randomly chosen, were supplied according to TDI-based
recommendations; ten of them were supplied as before) showed a
significant increase in gross yield: On average, the increase in the
gross yield in our blind study was around one percentage point.  The
probability that gross yield improvements of at least 0.25 percentage
points occurred is at least
87.6\,% according to the Wilcoxon rank sum test if inconsistent data is ignored
(even 95.5\,\% if inconsistent data is repaired in a plausible way).
And: this was the result of a single iteration of the optimization
procedure, which did not result in perfectly balanced Top-Dog-Indices.

Given the large economies of scale of a fashion discounter, we
consider the TDI a valuable contribution to revenue management tools
in this business sector.  Moreover, to the best of our knowledge, our
blind study is the first published study that evaluates a revenue
management method in the apparel retailer industry by comparing
simultaneously obtained business results of test-branches (optimized)
and control-branches (no action).

The draw-back of the TDI is its lack of information about the cardinal
expected revenue for a given size distribution of the supply.  This is
partly due to the fact that the loss of a bad size distribution is
closely related to the markdown policy of the discounter.  This
markdown policy, however, is itself subject of revenue management
methods.  Therefore, we regard the integration of size and price
optimization, as is done in our project BFS-DISPO, as a valuable
direction of further research .

\bibliographystyle{amsplain}
\bibliography{NKD-DO_ar}

\end{document}